\documentclass[11pt,reqno]{amsart}
\usepackage{amssymb,amscd,amsbsy,mathrsfs}
\usepackage{mathdots}
\renewcommand{\a}{\alpha}

\newcommand{\z}{\zeta}

\newcommand{\h}{{\mathscr H}}

\newcommand{\T}{{\Bbb T}}

\newcommand{\dd}{{\Bbb D}}
\newcommand{\R}{{\Bbb R}}
\newcommand{\Z}{{\Bbb Z}}

\newcommand{\0}{{\boldsymbol{0}}}

\newcommand{\bs}{\boldsymbol}

\newcommand{\m}{{\boldsymbol m}}
\newcommand{\bS}{{\boldsymbol S}}

\newcommand{\rf}[1]{(\ref{#1})}

\newcommand{\df}{\stackrel{\mathrm{def}}{=}}

\newcommand{\Ker}{\operatorname{Ker}}

\newcommand{\trace}{\operatorname{trace}}

\newcommand{\eeq}{\end{equation}}
\newcommand{\beq}{\begin{equation}}
\newcommand{\bay}{\begin{eqnarray}}
\newcommand{\ba}{\begin{align*}}
\newcommand{\ea}{\end{align*}}
\newcommand{\ey}{\end{eqnarray}}
\newcommand{\bey}{\begin{eqnarray*}}
\newcommand{\eey}{\end{eqnarray*}}

\newcommand{\be}{\infty}

\newcommand{\bl}{\blacksquare}

\newcommand{\im}{\operatorname{Im}}

\newtheorem{thm}{\hspace{\parindent}Theorem}[section]

\newtheorem{cor}[thm]{\hspace{\parindent}Corollary}
\newtheorem{lem}[thm]{\hspace{\parindent}Lemma}

\usepackage[T2A,T1]{fontenc}
\DeclareSymbolFont{cyrillic}{T2A}{cmr}{m}{it}
\def\makecyrsymbol#1#2{%
    \begingroup\edef\temp{\endgroup
        \noexpand\DeclareMathSymbol{\noexpand#1}
        {\noexpand\mathalpha}{cyrillic}%
        {\expandafter\expandafter\expandafter
            \calccyr\expandafter\meaning\csname T2A\string#2\endcsname\end}}%
    \temp}
\expandafter\def\expandafter\calccyr\string\char#1\end{#1}

\makecyrsymbol\Zhe\CYRZH
\makecyrsymbol\Be\CYRB
\makecyrsymbol\Shcha\CYRSHCH
\makecyrsymbol\Sha\CYRSH

\makeatletter
\def\upintkern@{\mkern-7mu\mathchoice{\mkern-3.5mu}{}{}{}}
\def\upintdots@{\mathchoice{\mkern-4mu\@cdots\mkern-4mu}%
 {{\cdotp}\mkern1.5mu{\cdotp}\mkern1.5mu{\cdotp}}%
 {{\cdotp}\mkern1mu{\cdotp}\mkern1mu{\cdotp}}%
 {{\cdotp}\mkern1mu{\cdotp}\mkern1mu{\cdotp}}}

\newcommand{\UpMultiIntegral}[1]{%
  \edef\ints@c{\noexpand\upintop
    \ifnum#1=\z@\noexpand\upintdots@\else\noexpand\upintkern@\fi
    \ifnum#1>\tw@\noexpand\upintop\noexpand\upintkern@\fi
    \ifnum#1>\thr@@\noexpand\upintop\noexpand\upintkern@\fi
    \noexpand\upintop
    \noexpand\ilimits@
  }%
  \futurelet\@let@token\ints@a
}
\makeatother

\DeclareFontFamily{OMX}{mdbch}{}
\DeclareFontShape{OMX}{mdbch}{m}{n}{ <->s * [0.8]  mdbchr7v }{}
\DeclareFontShape{OMX}{mdbch}{b}{n}{ <->s * [0.8]  mdbchb7v }{}
\DeclareFontShape{OMX}{mdbch}{bx}{n}{<->ssub * mdbch/b/n}{}

\DeclareSymbolFont{uplargesymbols}{OMX}{mdbch}{m}{n}
\SetSymbolFont{uplargesymbols}{bold}{OMX}{mdbch}{b}{n}
\DeclareMathSymbol{\upintop}{\mathop}{uplargesymbols}{82}
\DeclareMathSymbol{\upointop}{\mathop}{uplargesymbols}{"48}

\DeclareFontEncoding{MDB}{}{}
\DeclareFontFamily{MDB}{mdbch}{}
\DeclareFontShape{MDB}{mdbch}{m}{n}{ <->s * [0.8]  mdbchrmb }{}
\DeclareFontShape{MDB}{mdbch}{b}{n}{ <->s * [0.8]  mdbchbmb }{}
\DeclareFontShape{MDB}{mdbch}{bx}{n}{<->ssub * mdbch/b/n}{}
\DeclareFontSubstitution{MDB}{cmr}{m}{n}
\DeclareSymbolFont{mathdesignB}{MDB}{mdbch}{m}{n}%
\SetSymbolFont{mathdesignB}{bold}{MDB}{mdbch}{b}{n}%
\DeclareMathSymbol{\upintclockwise}{\mathop}{mathdesignB}{128}
\DeclareMathSymbol{\upointclockwise}{\mathop}{mathdesignB}{130}
\DeclareMathSymbol{\upointctrclockwise}{\mathop}{mathdesignB}{132}
\DeclareMathSymbol{\upoiint}{\mathop}{mathdesignB}{134}
\DeclareMathSymbol{\upoiiint}{\mathop}{mathdesignB}{136}

\makeatletter
\newcommand{\upint}{\DOTSI\upintop\ilimits@}
\newcommand{\upoint}{\DOTSI\upointop\ilimits@}
\makeatother

\theoremstyle{remark}

\newtheorem*{rem*}{Remark}

\newcommand\CA{{\rm C}_{\rm A}}

\newcommand{\ri}{{\rm i}}

\begin{document}

\numberwithin{equation}{section}

\numberwithin{equation}{section}

\title{Real-valued spectral shift functions for contractions and dissipative operators}
\author{M.M. Malamud, H. Neidhardt, V.V. Peller}
\thanks{The research on is supported by 
Russian Science Foundation [grant number 23-11-00153].}
\maketitle


\

\begin{abstract}
In recent joint papers the authors of this note solved a famous problem remained open for many years and proved that for arbitrary contractions with trace class difference there exists an integrable spectral shift function,
for which an analogue of the Lifshits--Krein trace formula holds. Similar results were also obtained for pairs of dissipative operators. Note that in contrast with the case of self-adjoint and unitary operators it may happen that there is no {\it real-valued} integrable spectral shift function. In this note we announce results that give sufficient conditions for the existence of an integrable real-valued spectral shift function in the case of pairs of 
contractions. We also consider the case of pairs of dissipative operators.
\end{abstract}

\setcounter{section}{0}
\section{\bf Introduction}
\setcounter{equation}{0}
\label{In}

\

In this note we study conditions, under which a pair of contractions on Hilbert space with trace class difference has a real-valued integrable spectral shift function on the unit circle $\T$. We also consider the case of pairs of dissipative operators.

Recall that  physicist I.M. Lifshits 
for a pair of self-adjoint operators 
$\{A_0,A_1\}$ in Hilbert space with trace class difference, introduced in \cite{L} a spectral shift function
$\bs\xi=\bs{\xi}_{A_0,A_1}$ on the real line $\R$ and discovered the trace formula
\bay
\label{foslsso}
\trace\big(f(A_1)-f(A_0)\big)=\int_\R f'(t)\bs{\xi}_{A_0,A_1}(t)\,dt
\ey
for sufficiently nice functions $f$ on $\R$.

Later M.G. Krein gave a rigorous mathematical justification of this formula and showed that in the case when 
$A_1-A_0$ is a trace class operator, the spectral shift function must be integrable on $\R$ and trace formula \rf{foslsso} holds for sufficiently nice functions $f$.


On the other hand, Krein observed in \cite{Kr} that the right-hand side of \rf{foslsso} is well-defined for an arbitrary Lipschitz function $f$ and posed the problem to describe the maximal class of functions  $f$, for which trace formula \rf{foslsso} holds for arbitrary self-adjoint operators with trace class difference. In particular, Krein asked whether it is possible to extend trace formula
\rf{foslsso} to the class of Lipschitz functions.

It turned out, however, that this is not the case. Indeed, Yu.B. Farforovskaya showed that there exist a Lipschitz function $f$ and self-adjoint operators $A_0$ and $A_1$ such that $A_1-A_0\in\bS_1$ but
$f(A_1)-f(A_0)\not\in\bS_1$\footnote{we denote by $\bS_1$ the trace class, see \cite{GK}}.


The Krein problem was completely solved in the paper \cite{Pe6}, in which it was shown that the maximal class of functions $f$, for which trace formula \rf{foslsso} holds for arbitrary self-adjoint operators with trace class difference coincides with the class of {\it operator Lipschitz functions} 
(see the survey \cite{AP}, which contains a lot of information on operator Lipschitz functions).

In \cite{Kr2} for pairs of unitary operators with trace class difference Krein constructed a spectral shift function
that is integrable and real-valued on $\T$ and established an analogue of trace formula \rf{foslsso}.

We proceed now to a trace formula for functions of contractions with trace class difference. 
Recall that an operator $T$ on Hilbert space is called a {\it contraction} if $\|T\|\le1$.

Recall also that the Sz.-Nagy--Foia\c s functional calculus (see \cite{SNF}) for a contraction $T$ on Hilbert space, associates with each function $f$ in the {\it disc-algebra} 
$\CA$ a function $f(T)$ of $T$; moreover, this functional calculus
$$
f\mapsto f(T),\quad f\in\CA,
$$
is linear and multiplicative. Also, the von Neumann inequality
$$
\|f(T)\|\le\|f\|_{\CA}\df\max\{|f(\z)|:~|\z|\le1\},\quad f\in\CA.
$$
holds.

A lot of attempts has been undertaken to generalize the Lifshits--Krein trace formula to the case of functions of contractions. The problem is for a pair of contractios $\{T_0,T_1\}$ with trace clas difference, to prove the existence of an integrable function $\bs\xi$  on the circle $\T$, for which the following trace formula would hold: 
\bay
\label{forsleszha}
\trace\big(f(T_1)-f(T_0)\big)=\int_\T f'(\z)\bs{\xi}(\z)\,{\rm d}\z
\ey
for sufficiently nice functions $f$ (for example, for analytic polynomials $f$). It would be natural to call such a function $\bs\xi$ a {\it spectral shift function}for the pair of contractions $\{T_0,T_1\}$. It is worth mentioning that if such a function $\bs\xi$ exists, it is by no means unique. 
Indeed, we can add to such a function $\bs\xi$ an arbitrary function in the Hardy class $H^1$ and obtain a new spectral shift function.

Such problems are also important for pairs of maximal dissipative operators. Recall that a densely defined operator $L$ is called {\it disssipative} if $\im(Lx,x)\ge0$ for any vector $x$ in the domain ${\rm D}(L)$ of $L$. A dissipative operator $L$ is called {\it maximal}, if it has no proper dissipative extension.

To describe the history of the problem, we start with Langer's paper \cite{La} of 1965, whose results imply the existence of an integrable spectral shift function under the assumption that the spectra of the contractions are contained in the open unit disc $\dd$. We also mention Rybkin's papers \cite{Ryb84,Ryb94}, the paper by Adamyan and Neidhardt \cite{AN} and also Krein's paper \cite{Kr87}. Note that in \cite{Ryb94} under additional assumptions on a pair of contractions, Rybkin proved the existence of an A-integrable complex-valued spectral shift function that is not necessarily Lebesgue integrable.

We underline here that in the Adamyan--Neidhardt paper \cite{AN} the existence of a real integrable spectral shift function for  pair of contractions $\{T_0,T_1\}$ was proved under a stronger assumption than $T_1-T_0\in\bS_1$; namely, under the assumption 
\bay
\label{logarifmy}
\sum_{k\ge0}s_k(T_1-T_0)\log\big(1+(s_k(T_1-T_0))^{-1}\big)<\be
\ey
(it is assumed that the function $x\mapsto x\log(1+x^{-1})$ takes value 0 at $x=0)$. It follows easily from our Theorem \ref{veshch} in \S\;2  that condition \rf{logarifmy} is not necessary for the existence of a real integrable spectral shift function.

The problem of obtaining an analogue of the Lifshits--Krein trace formula for functions of contractions remained open for many years and was completely resolved in  \cite{MNP1} (see also the paper 
\cite{MN}, in which the trace formula was obtained under an additional assumption). Another solution to this problem was obtained in \cite{MNP2}. 
Moreover, in \cite{MNP1} and \cite{MNP2} the authors described the maximal class of functions $f$, for which trace formula \rf{forsleszha} holds for arbitrary pairs of contractions with trace class difference. This class coincides with the class of operator Lipschitz functions analytic in the disc $\dd$.

It is known that (see \cite{MNP1}) that for contractions with trace class difference, it is not always true that there exists an integrable real-valued spectral shift function. Nevertheless, as shown in \cite{MNP2}, a pair of contractions with trace class difference always has a real-valued A-integrable spectral shift function. 

In this note we announce a sufficient condition for a pair of contractions with trace class difference possesses
a real-valued integrable spectral shift function. Recently, in the paper \cite{CS} by Chattopadhyay and Sinha
it was shown that if  $T_0$ and $T_1$ are contractions with trace class difference and $T_0$ is a strong contraction, i.e., $\|T_0\|<1$, then the pair $\{T_0,T_1\}$ has an integrable real-valued spectral shift function. In
this note by developing ideas in the paper \cite{MNP2} we essentially improve the result of \cite{CS}.
This will be done in \S\:\ref{szhali}.

It is also worth mentioning here that in \cite{MNP2} it was shown that if $T$ is a contraction and $U$ is a unitary operator such that $U-I\in\bS_1$, then the pair $\{T,UT\}$ has a real integrable spectral shift function (see Lemma 9.1); on the other hand, 
if $T$ is a contraction and $X$ is a contraction such that 
$X\ge\0$ and $I-X\in\bS_1$, then the pair $\{T,XT\}$ has a purely imaginary integrable spectral shift function (see Lemma 9.2).

The Lifshits--Krein trace formula was also generalized in \cite{MNP1} and \cite{MNP2} to the case of maximal dissipative operators. If
$L_0$ and $L_1$ are maximal dissipative operators with trace class difference, then as shown in \cite{MNP1} and \cite{MNP2}, there exists an integrable spectral shift function $\bs\xi$ on the real line $\R$ such that the trace formula
\bay
\label{fcddo}
\trace\big(f(L_1)-f(L_0)\big)=\int_\R f'(t)\bs{\xi}(t)\,{\rm d}t
\ey
holds at least for rational functions $f$ with poles in the open lower half-plane.

On the other hand, if we replace the condition $L_1-L_0\in\bS_1$ with the resolvent condition
$$
(\ri I+L_1)^{-1}-(\ri I+L_0)^{-1}\in\bS_1,
$$
then there exists a spectral shift function $\bs\xi$ on $\R$ such that
\bay
\label{s_vesom}
\int_\R|\bs{\xi}(t)|(1+t^2)^{-1}\,{\rm d}t<\be
\ey
and the trace formula \rf{fcddo} holds for rational functions $f$ with poles in the open lower half-plane
(see \cite{MNP2}). 

In \S\:\ref{dissipatsiya} of this note we study the problem of the existence of a real-valued spectral shift function for pairs of dissipative operators. 

\medskip

\section{\bf Real-valued integrable spectral shift functions for contractions}
\setcounter{equation}{0}
\label{szhali}

\medskip

To state a sufficient condition for the existence of a real-valued integrable spectral shift function, we need the notion 
of the defect operators of contractions. If $T$ is a contraction, its defect operators are defined by 
$$
D_T \df(I - T^*T)^{1/2},\quad\mbox{and}\quad D_{T^*}\df(I - TT^*)^{1/2}.
$$

To obtain the main result of this section, we need the following lemma, in which $\bS_p$ denotes the Schatten--von Neimann class, see \cite{GK}.

\begin{lem}
\label{raz def}
Let $0<p<\be$. Suppose that
$T_0$ and $T_1$ are contractions on Hilbert space such that
\bay
\label{yadra}
\Ker D_{T_0} = \{\0\},
\ey
and
\begin{equation}
\label{-2a}
(T_1-T_0)D_{T_0}^{-2\a}\in\bS_p
\quad \mbox{and} \quad
(T_1^*-T^*_0)D_{T^*_0}^{-2\a}\in \bS_p
\end{equation}
for some $\a$ in $(\tfrac{1}{2},1]$.
Then
\bay
\label{defekty}
D_{T_1}-D_{T_0}\in \bS_p\quad\mbox{and}\quad D_{T_1^*} - D_{T^*_0} \in \bS_p.
\ey
\end{lem}

It is easy to see that each of the inclusions in \rf{-2a} implies that $T_1- T_0 \in\bS_p$.
Note also that 
$\Ker D_{T_0} = \{\0\}$ if and only if $\Ker D_{T^*_0} = \{\0\}$, and so, condition \rf{yadra} can be replaced with the condition $\Ker D_{T^*_0} = \{\0\}$.


The following assertion is the main result of this section. 

\begin{thm}
\label{veshch}
Suppose that $T_1$ and $T_0$ are contractions satisfying {\em{\rf{yadra}}} and
{\em{\rf{-2a}}} with $p=1$. Then the pair $\{T_0,T_1\}$ possesses a real-valued spectral shift function.
\end{thm}

{\bf The idea of the proof.} To deduce Theorem \ref{veshch} from Lemma \ref{raz def}, we are going to use the Schaffer matrix  unitary dilation. If $T$ is a contraction on a Hilbert space $\h$, consider the unitary operator   $U^{[T]}$ on the space $\ell^2_\Z(\h) = \bigoplus_{j \in \Z}\h_j$,  $\h_j = \h$, of two-sided $\h$-valued sequences, see
\cite{SNF}, Ch. 1, \S\:5 (see also the paper \cite{ MNP2}, in which the Schaffer dilations were used in connection with spectral shift functions for contractions). The space $\h$ is identified with the subspace of the space of sequences of the form $\{v_n\}_{n\in\Z}$ such that $v_j=\0$ with $j\ne0$. 
The operator $U^{[T]}$ is defined by the operator matrix
\bay 
\label{maShe}
U^{[T]}=
\left(\begin{matrix}
\ddots&\ddots&\vdots&\vdots&\vdots&\vdots&\vdots&\vdots&\iddots\\
\cdots&\0&I&\0&\0&\0&\0&\0&\cdots\\
\cdots&\0&\0&I&\0&\0&\0&\0&\cdots\\
\cdots&\0&\0&\0&D_T&-T^*&\0&\0&\cdots\\
\cdots&\0&\0&\0&T&D_{T^*}&\0&\0&\cdots\\
\cdots&\0&\0&\0&\0&\0&I&\0&\cdots\\
\cdots&\0&\0&\0&\0&\0&\0&I&\cdots\\
\iddots
&\vdots&\vdots&\vdots&\vdots&\vdots&\vdots&\ddots&\ddots
\end{matrix}\right).
\ey
In this matrix the entry $T$ has position $(0,0)$. 
In other words, the entries $U^{[T]}_{j,k}$ of the operator matrix
$U^{[T]}$ are defined by
$$
U^{[T]}_{0,0}=T,\quad U^{[T]}_{0,1}=D_{T^*},\quad U^{[T]}_{-1,0}=D_T,
\quad U^{[T]}_{-1,1}=-T^*,\quad
U^{[T]}_{j,j+1}=I\quad\mbox{for}\quad j\ne0,~-1,
$$
while the remaining entries are equal to $\0$.

It is easy to see that $U^{[T]}$ is a unitary dilation of $T$.
Clearly, if $T_0$ and $T_1$ satisfy the hypotheses of Lemma \ref{raz def} with $p=1$, then the unitary operators  $U^{[T_0]}$ and $U^{[T_1]}$ have trace class difference, and so they have a real spectral shift function. It is also easy to see that a spectral shift function for the pair $\big\{U^{[T_0]},U^{[T_1]}\big\}$ is also a spectral shift function  for the initial  pair $\{T_0,T_1\}$ which implies the conclusion of Theorem \ref{veshch}. $\bl$

\medskip

Recall that by the Sz.-Nagy--Foia\c s theorem (see \cite{SNF}), the spectral measure of the minimal unitary dialtion of a completely nonunitary contraction must be mutually absolutely continuous  with respect to Lebesgue measure on the circle.
Though the Schaffer matrix dilation does not have to be minimal, its spectral measure still must be
mutually absolutely continuous  with respect to Lebesgue measure (see \cite{MNP2}, Corollary 10.2).

\medskip

{\bf Remark.} It is obvious that in the case $\|T_0\|<1$ the pair $\{T_0,T_1\}$ satisfies the hypotheses  of Theorem \ref{veshch}, and so Theorem \ref{veshch} improves the result of \cite{CS} that has been mentioned in the introduction.

\begin{cor}
Let $T$ and $X$ be contractions such that $X\ge\0$ and $I-X\in\bS_1$. Suppose that  
$$
\Ker D_T = \{\0\},
$$
and
$$
(I-X)TD_T^{-2\a}\in\bS_1
\quad \mbox{and} \quad
T^*(I-X)D_{T^*}^{-2\a}\in \bS_1
$$
for some $\a$ in $(\tfrac{1}{2},1]$.
Then the pair $\{T,XT\}$ has both a real integrable spectral shift function and a purely imaginary integrable 
spectral shift function.
\end{cor}

Indeed, the first part of the conclusion of the corollary is an immediate consequence of Theorem \ref{veshch} while the second part is just Lemma 9.2 of \cite{MNP2} that has been mentioned in the introduction.

Note also that if $\bs\xi$ is a real-valued spectral shift function and $\ri\bs\eta$ is a purely imaginary spectral shift function, then  the function $\bs\xi-\ri\bs\eta$ belongs to the Hardy class $H^1$. In this case the harmonically conjugate functions $\widetilde{\bs\xi}$ and $\widetilde{\bs\eta}$ are integrable. If, in addition to this, $\bs\xi\ge\0$, then the function $\bs\xi$ satisfies the Zygmund condition
$$
\int_\T\bs\xi(\z)\log(1+\bs\xi(\z))\,{\rm d}\m(\z)<\be.
$$
Recall in this connection that in \cite{MN} and \cite{MNP1} it was observed that if a pair of contractions with trace class difference has a complex-valued spectral shift function that satisfies the Zygmund condition, then there is also a real-valued spectral shift function. 

To conclude the section, we would like to mention that the proof of Lemma \ref{raz def} is based on the following assertion.

\begin{lem}
\label{X,Y}
Let $0<p<\be$ and $\a\in(\tfrac{1}{2},1]$. 
Suppose that $X$ and $Y$ are operators on Hilbert space such that $\0 \le X\le I$, $\0 \le Y\le I$ and 
$\Ker X = \{\0\}$. If $Y-X\in\bS_p$ and the operator
$(Y-X)X^{-\a}$ extends to the operator of class $\bS_p$, then $Y^{1/2} - X^{1/2} \in \bS_p$.
\end{lem}

\medskip

\section{\bf The case of dissipative operators}
\label{dissipatsiya}

\medskip

In this section we give a sufficient condition on a pair $\{L_0,L_1\}$ of maximal dissipative operators with trace class difference for the existence of a real-valued spectral shift function. Nevertheless, our condition leads to the existence of a spectral shift function $\bs\xi$ satisfying \rf{s_vesom}. It turns out that under our assumptions 
a real-valued {\it integrable} function does not have to exist. Moreover, such a function cannot exist if
$\trace(L_1-L_0)\not\in\R$.

Throughout this section we use the notation 
$$
V \df L_1-L_0.
$$

\begin{thm}
\label{imL0}
Let $L_0$ and $L_1$ be bounded dissipative operators such that $V = L_1-L_0\in\bS_1$ and
$\Ker\im L_0=\{\0\}$.  If
\bay
\label{Im-1}
(L_1- L_0)(\im L_0)^{-1} \in \bS_1,
\ey
then the pair $\{L_0,L_1\}$ possesses a real-valued spectral shift function $\bs\xi$ that satisfies the condition 
\bay
\label{svesom}
\int_\R|\bs{\xi}(t)|(1+t^2)^{-1}\,{\rm d}t<\be.
\ey
\end{thm}

Note that by Theorem 5.6 of \cite{MN}  
for bounded dissipative operators $L_0$ and $L_1$ under the assumption
$V\in\bS_1$ it is always true (i.e. without condition (3.1)) that there exists an {\it integrable} complex-valued  spectral shift function $\bs\xi$. Moreover, it was shown in 
 \cite{MNP2} (see Theorem 9.6 and Corollary 9.7) that an integrable  complex-valued spectral shift function  exists for pairs of dissipative with trace class difference even without the assumption of their boundedness. 
Also, if $V$ is dissipative, there exists a spectral shift function  
$\bs\xi$ satisfying the inequality $\im \bs\xi\ge \0$. If $V$ is self-adjoint, a spectral shift function $\bs\xi$ can be chosen to be real-valued while in the case $V = -V^*$ there exists a purely imaginary spectral shift function.

Theorem \ref{imL0} complements the above result by stating that under condition \eqref{Im-1} 
it is always true (i.e., even without the assumption that $V$ is self-adjoint) that it is possible to make a spectral shift function  {\it real-valued}  and satisfying  \rf{svesom} but not necessarily integrable.  In particular, in the case $V = -V^*$, for a pair  $\{L_0,L_1\}$ under condition \eqref{Im-1} there exists a purely imaginary integrable
spectral shift function  $\bs\xi_1$ and also a real-valued not necessarily integrable spectral shift function  
$\bs\xi_2$ satisfying \rf{svesom}.

As one can see from the following theorem, it is very far from being true that an integrable real-valued spectral shift function always exists.

\begin{thm}
Suppose that $\{L_0,L_1\}$ are pairs of dissipative operators satisfying the hypotheses of Theorem {\em\ref{imL0}}. Suppose that $\trace(L_1-L_0)\not\in\R$.
Then the pair $\{L_0,L_1\}$ cannot have an integrable real-valued spectral shift function.

Moreover, if $V$  is dissipative, then the condition $\trace(L_1-L_0) \in\R$ is equivalent to the existence
of an integrable real-valued spectral shift function.
\end{thm}

\
 
 \begin{footnotesize}
 
\noindent
\begin{tabular}{p{8cm}p{15cm}}
M.M. Malamud & V.V. Peller \\
People's Friendship University & St.Petersburg State University\\
of Russia (RUDN University) & Universitetskaya nab., 7/9\\
6 Miklukho-Maklaya St, Moscow & 199034 St.Petersburg\\
Russsia&Russia\\
email: malamud3m@gmail.com 
\\

&St.Petersburg Department\\
&Steklov Institute of Mathematics \\
&Russian Academy of Sciences\\
&Fontanka 27, 191023, 191023 St.Petersburg\\
&Russsia\\
& email: peller@math.msu.edu
\end{tabular}

\end{footnotesize}

\

\end{document}